\documentclass{amsart}[14pt]

\usepackage{amsmath, amstext, amssymb, amsopn, amsthm, amscd, amsfonts}

\def\dobleint{\mathop{\int \mkern-8mu \int}}

\newcommand{\C}{\mathbb C}

\newcommand{\D}{\mathbb D}

\newcommand{\N}{\mathbb N}

\newcommand{\chat}{\widehat{\mathbb C}}

\newtheorem{theo}{Theorem}

\newtheorem{lema}{Lemma}

\newtheorem{prop}{Proposition}

\newtheorem{corol}{Corollary}

\newtheorem{defi}{Definition}

\newtheorem{rema}{Remark}

\begin{document}

\title{\bf Poincare Series and instability of exponential maps}

\author{P.Makienko and G. Sienra}

\maketitle

\begin{abstract}
We relate the properties of the postsingular set for the exponential family to the
questions of stability. We calculate the action of the Ruelle operator for the
exponential family. We prove that if the asymptotic value is a summable point and
its orbit satisfies certain topological conditions, the map is unstable hence there
are no Beltrami differentials in the Julia set. Also we show that if the
postsingular set is a compact set, then the singular value  is summable.
\end{abstract}

%%\keywords{exponential map, stability, summability, Ruelle operator.}

\bigskip

\section{Introduction}

If $f$ is a transcendental entire map, we denote by $f^{n}$, $n \in {\N}$, the n-th
iterate of $f$ and write the Fatou set as $F(f)=\{z \in {\C}$; there is some open
set U containing $z$ in which $\{f^{n}\}$ is a normal family \}. The complement of
$F(f)$ is called the Julia set $J(f)$. We say that $f$ belongs to the class $S_{q}$
if the set of singularities of $f^{-1}$ contains at most $q$ points.

Two entire maps $g$ and $h$ are topologically equivalent if there exist
homeomorphisms $\varphi$, ${\psi}:{\C} \rightarrow {\C}$ such that ${\varphi} \circ
g = h \circ {\psi}$. Given a map $f$, let us denote by $M_{f}$, the set of all
entire maps topologically equivalent to $f$.

 It is proved in \cite{EL} that $M_{f}$ has the structure of a (q+2)-dimensional complex manifold.
The Affine group acts on the space $M_{f}$ and as it shown in \cite{EL} the space $
N_f = M_{f}/\{{\text{\sl Affine group}}\} $ is a $ q-$dimensional complex orbifold.

A measurable field of tangent ellipses of bounded eccentricity determines a complex
structure on the sphere. This ellipse field is recordered by a (-1,1)-form
${\mu}(z)\frac{d{\overline{z}}}{dz}$ with $\vert \vert {\mu}\vert\vert_{\infty}<1$,
a {\sl Beltrami differential}. If an entire map $ f $ is holomorphic in a complex
structure defined by the Beltrami differential $ \mu, $ then $ \mu $ is the {\sl
invariant Beltrami differential}. Since the sphere admits a unique complex
structure, there is a homeomorphism ${\phi}:{\C}\rightarrow{\C}$ such that ${\mu}$
is the pullback of the standard structure and the map $ f_{\phi} = \phi\circ
f\circ\phi^{-1} $ is an entire map.
%% through the Beltrami operator.

The non existence of an invariant Beltrami differential (invariant line field) on
the Julia set is related to the Fatou conjecture, see \cite{MM}.
%%The approach to the
%%non existence of line fields in [5] and [6] is for rational maps and in [3] for
%%entire transcendental maps with algebraic  singularities.

% by considering the action of the Ruelle

%operator.

Now let us consider the main hero of this paper - {\sl Exponential family}: $ E =
\{f_{\lambda}(z) = \exp(\lambda z), \lambda \in {\C^*}\}. $ Then  $ N_{f_1} \cong E,
$ where  $  f_1 =\exp(z). $ The map $ f_{\lambda_0} $ is {\sl structurally stable}
if for any $ \lambda $ close enough to $ \lambda_0 $ there exists a quasiconformal
homeomorphism $ \phi_\lambda, $ such that $ f_\lambda = \phi_\lambda\circ
f_{\lambda_0}\circ \phi_\lambda^{-1}. $

Due to Man{\'e}, P. Sad, D. Sullivan (see \cite{MSS}) and A. Eremenko, M. Lyubich
(see \cite{EL}) the following three items are equivalent for $ E $:

\begin{itemize}
\item  Fatou conjecture \item There is no invariant Beltrami differentials supported
by the Julia set \item If $ J(f_{\lambda}) = {\C}, $  then $ f_{\lambda} $ is
structurally unstable.
\end{itemize}

In 1985 R. Devaney (see \cite{D}) proves that $ \exp(z) $ is structurally unstable,
after A. Douady and L. R. Goldberg (see \cite{DG}) did show that the maps $
\lambda\exp(z), \lambda \geq 1 $ are topologically unstable. Zhuan Ye (see
\cite{Ye}) proves that $ f_{\lambda} $ is structurally unstable map if $
\lim_{n\to\infty} f_{\lambda}^n(0) = \infty. $

In this paper we follow  the approach of papers \cite{Av}, \cite{Lev} and
\cite{Mak}-\cite{Mak1}, (case of rational maps) and \cite{DMS} (case of
transcendental entire maps with only algebraic singularities). In the case of
Exponential family we have only one asymptotic singularity which is a different
situation that in \cite{DMS}.

The stability of a map depends on the behavior of the {\it postsingular set},
denoted as $X_{\lambda} =\overline{\{\cup_{n\geq 1} f_\lambda^n(0)\}}. $
%%Therefore
%%our first theorems relate $ z = 0 $ and its orbit in the set $X_{\lambda}$ with the
%%stability or instability of the map in question.

Let us start with $ f_{\lambda} $ whose Julia set is equal to the plane. Then  we
have the following simple possibilities:
\begin{enumerate} \item $lim_{n\to\infty}
\vert (f_\lambda^n)'(0)\vert = 0, $ \item there exists a subsequence $ \{n_i\} $
such that $ lim_{i\to\infty} \vert (f_\lambda^{n_i})'(0)\vert = \infty, $ \item
there exists a subsequence $ \{n_i\} $ such that $ lim_{i\to\infty} \vert
(f_\lambda^{n_i})'(0)\vert = M < \infty $ and $ M \neq 0.$
\end{enumerate}

We believe that the first case contains  a contradiction. Since in this situation
the forward orbit of $ 0 $ must converge to an attractive cycle and hence $ 0 \notin
J(f_{\lambda}).$ We show this conjecture under very strong additional conditions
only as an illustration that this conjecture is not completely false (see theorem
1).

As for the last two cases, the Fatou conjecture claims that  $ f_\lambda $ is an
unstable map.
%%The logical order of our theorems is as follows:

Define
\begin{defi} Let $ \lambda \in {\C}^*, $ then the {\bf Poincar{\'e} series} for $ f_\lambda $
is the following formal series
$$
P_\lambda =1 + \frac{1}{\lambda}\sum_{i = 2}^{\infty}\frac{1}{(f_{\lambda}^{i -
2})^\prime(1)}.
$$
\end{defi}

Let
$$
S_{n}=1+\frac{1}{\lambda}\sum_{i=2}^{n}\frac{1}{(f_{\lambda}^{i-2})^\prime(1)}
$$
be a particular sums of the Poincar{\'e} series $ P_\lambda. $ Then we have the
following theorem
\begin{theo}\hfill
\begin{enumerate}
\item If there exist a sequence $ \{n_{i}\} $ such that $
(f_{\lambda}^{n_i})^\prime(1) \to \infty $ and $ lim_{i \to\infty}sup\vert
S_{n_i}\vert>0, $ then $ F(f_{\lambda}) =\emptyset$ and $f_{\lambda}$ is unstable.
\item If there exist a sequence $\{n_{i}\}$ such that
$(f_{\lambda}^{n_i})^\prime(1)\asymp c$, where $c\neq 0$ is a constant and $lim_{i
\to\infty}sup|S_{n_{i}}|=\infty$, then $f_{\lambda}$ is unstable. \item Let $
\lim_{n\to\infty} (f^n_{\lambda})^\prime(1) =0, $ and one of the following
conditions holds:

$\lim_{n\to \infty}sup\frac{\vert (f^{n +1}_{\lambda})^\prime(1)\vert}{\vert
(f^n_{\lambda})^\prime(1)\vert} < \infty, $ or

$\lim_{n\to \infty}inf\frac{\vert (f^{n +1}_{\lambda})^\prime(1)\vert}{\vert
(f^n_{\lambda})^\prime(1)\vert} > 0. $

Then $ F(f_{\lambda}) \neq \emptyset. $
\end{enumerate}
\end{theo}
\begin{prop}Does not exist a map $f_{\lambda}$ such that $lim_{n
\rightarrow \infty}|f_{\lambda}^{n'}(1)|=C>0$
\end{prop}
The next theorems discuss the best conditions on the Poincar{\'e} series  and on the
postsingular set for the map to be unstable.

%%Let $X_{\lambda}= \overline{\bigcup_{n=0}^{\infty}f_{\lambda}^{n}(0)}$.

\begin{defi}
A  point $a \in \C$ is called "{\rm summable}" if and only if the series
$$
\sum_{i = 0}\frac{1}{(f_{\lambda}^i)'(a)}
$$
is absolutely convergent. Note that the point $ z = 0 $ is summable if and only if
the Poincar{\'e} series $ P_{\lambda} $ is absolutely convergent.
%We say that
%$f_{\lambda}$ is summable if $ a = 1 $ is a summable point.
\end{defi}

\begin{defi} Let $ W \subset E $ be the subset of exponential maps $ f_{\lambda}, $
with summable singular point $ 0 \in J(f_{\lambda}), $ satisfying one of the
following conditions:
\begin{enumerate}
\item $ 0 \not\in X_{\lambda}, $ \item $ X_{\lambda} $ does not separate the plane,
\item $ m(X_{\lambda}) = 0, $ where $ m $ is the Lebesgue  measure.
\end{enumerate}
\end{defi}

\begin{theo} Let $ f_{\lambda} \in W $. Then $ f_{\lambda} $ is an unstable map, and hence there is
no invariant Beltrami differentials on its Julia set.

\end{theo}
\begin{theo} Let $ f_{\lambda} \in E, $ with $ J(f_{\lambda}) = {\C}. $
Assume $ 0\notin X_{\lambda}$ (i.e. $0$ is non-recurrent), then
\begin{enumerate}
\item There exist a subsequence $n_{k}$ such that $(f^{n_k})^\prime(1) \rightarrow
\infty$ \item If $ X_{\lambda} $ is bounded, then the singular point $ z = 0 $  is
summable for $ f_{\lambda}. $
\end{enumerate}
\end{theo}
 In section 2 we discuss and prove Theorem 3 and Proposition 1.

 In section 3 we consider the basic definitions and properties of the Ruelle operator and
the potential of deformations, as a consequence we prove theorem 1.

 The rest of the paper is devoted to prove Theorem 2.
 %%Theorems related to the hyperbolic, hence stable behavior of this family have been
 %%studied in [1]and in [2] Devaney studies the unstable behavior
%%of the exponential maps with $\lambda \in {\R}$.

\section{Postsingular set and dynamics}

Ma\~{n}e has a result that establishes expansion properties of  rational maps on the
compact subsets of their Julia sets, which are far away from the parabolic points
and the $w$-limit sets of recurrent critical points. Next we will consider this
result for our map $f_{\lambda}$.

\begin{rema}
Note that if $f_{\lambda}^{n}(0)\rightarrow \infty$ then $f_{\lambda}$ is summable.
To see this, consider
$$
|{\frac{1}{(f_{\lambda}^{n+1})^\prime(a)}}|/|\frac{1}{(f_{\lambda}^{n})^\prime(a)}|
= |\frac{1}{\lambda f_{\lambda}^{n}(a)}|
$$
now choose $ a = f_{\lambda}(1)$ and since the orbit of $0$ tends to $\infty$ this
fraction converges to zero, so the series $ \sum
\frac{1}{(f_{\lambda}^{n})^\prime(1)}$ absolutely converges.
\end{rema}
\bigskip
\subsection{Proof of Theorem 3}The proof of the theorem follows exactly the proof
in \cite{ST}, by Shishikura and Tan Lei. For completeness we will state the lemmas
used in the paper above mentioned, restricted to the situation of our case. Hence in
order to prove our theorem 3, we will follow their arguments.

Denote by $ d(z,E) $ the Euclidian distance between a point $z\in \C$ and a closed
subset $ E \subset \C. $ Let $ d_Y(z, X) $ be the Poincar{\'e}  distance on a
hyperbolic surface $ Y $ between a point $ z $ and a closed subset $ X \subset Y $
and $ diam_{W}(W') $ the diameter of $ W' $ with respect of the the Poincar{\'e}
metric of $ W. $

\begin{lema}(\cite{ST} lemma 2.1). For any $ 0 \leq r \leq 1, $ there exist a constant
$ C(1,r) \geq 0 $ such that for any holomorphic proper map $ g : V \rightarrow \D $
of degree 1, with $ V $ simple connected, each component  of $
g^{-1}(\overline{D_{r}(0)})$ has diameter $ \leq C(1,r) $ with respect to the
Poincar{\'e} metric on $ V. $ Moreover $ lim_{r \rightarrow 0}C(1,r) = 0. $
\end{lema}

\begin{defi} $ {\bf N_{0}} $: There exist $z_{1},...,z_{N_{0}-1} \in \D $
such that $ \{\frac{2}{3}\leq |z| \leq 1\} \subset \bigcup_{i=1}^{N_{0}-1}
D_{\frac{1}{3}}(z_{i}). $ Let $ C_{0}=N_{0}C(1,\frac{2}{3}) $
\end{defi}

The Julia set  $ J(f_{\lambda}) = \C, $ hence we can choose  a periodic point $ w $
so that the domain $ \Omega = \C\backslash \{\text{forward orbit of the point w}\} $
satisfies: $ d_{\Omega}(0,X_{\lambda}) \geq 2C_{0}. $

\begin{lema}(\cite{ST} lemma 2.3) Let $ U_{0}=D_{r}(x) $ be a disc centered at $ x \in X_{\lambda} $
with radius $ r $ so that   $ U_{0} \subset \Omega $ and $ diam_{\Omega} (U_{0})
\leq C_{0}, $ then for every $n \geq 0 $ the following is true:

${\bf deg(n)}$. For every $ D_{s}(z) \subset U_{0} $ with $ 0 \leq s \leq
d(z,{\partial}U_{0})/2, $ and every connected component $ V' $ of $
f_{\lambda}^{-n}(D_{s}(z)), $ $ V' $ is simply connected and $
deg(f_{\lambda}^{n}:V' \rightarrow D_{s}(z))=1; $

${\bf diam(n)}$. For every $ D_{r}(w) \subset U_{0} $ with $ 0 \leq r \leq d(w,
{\partial} U_{0})/2 $ and every connected component of $ V $ of $ f_{\lambda}^{-n}
(D_{r}(w)), $ $diam_{\Omega}V \leq C_{0}. $

\end{lema}

Now we begin to prove the theorem 3.  If only $ \infty $ is a point of accumulation
of $ \{\cup_n f_{\lambda}(0)\}, $ then by the remark 1 above the point $ z = 0 $ is
a summable and hence $ \lim_{n\to\infty}\vert(f_{\lambda}^n)^\prime(0)\vert =\infty.
$

Now let $ y \in X_{\lambda} $ be another point of accumulation of the orbit of $ z =
0. $ Let $ n_i $ be any subsequence such that $ y = \lim_{i\to\infty}
f^{n_i}_{\lambda}(1). $ Then we claim:

{\bf Claim} $ \lim_{i\to\infty}\vert (f_{\lambda}^{n_i})^\prime(1)\vert = \infty. $

To prove the claim we repeat the arguments of Shishikura and Tan Lei. Assume there
exist a number $ M < \infty $ and a sequence of natural numbers $\{n_{j}\} \subset
\{n_i\} $ such that $ |(f_{\lambda}^{n_j})^\prime(1)| \leq M. $ Then by the lemma 2
there exist an integer $ N $ and a number $ r $ such that components $ W_j \subset
f_{\lambda}^{-n_j} (D_r(y)) $ containing the point $ z =1 $ are simply connected and
the respective restriction maps $ f_{\lambda}^{n_j}: W_J \rightarrow D_{r}(y) $ are
univalent for all $ j \geq N. $ Now let $ B \subset \Omega $ be the hyperbolic ball
of the radius $ C_0 $ centered at the point $ z = 1, $  then $ B $ is a precompact
subset of $ \Omega $ and hence has a bounded Euclidian diameter in $ {\C.} $
Besides, again by the lemma 2, the set $ \{\cup_j W_j\} \subset B. $ Let $ g_{j}:D
\rightarrow W_{j} $ be the inverse maps, then it is a normal family. Hence after
passing to a subsequence we cam assume that $ g_j $ converge. Let  $ g_{\infty} $ be
a limit map, then $ g_{\infty} \neq const $ since the derivatives are $\geq
\frac{1}{M}$ by hypothesis. Then there is a neighborhood $ U_{0} $ of $ z = 1 $ such
that $ U_{0} \subset g_{j}(D) $ for large $ {j}. $ Then $ f_{\lambda}^{n_{j}} $ is a
normal in $ U_0, $ but there are many periodic expansive points in $ U_0 \subset
J(f_\lambda) $ and the derivative diverges. Which is a contradiction. The claim and
the first part of the theorem are done.

Finally for the proof of the second part, we again repeat arguments of Shishikura
and Tan Lei in \cite{ST}. So assume that $ f_{\lambda} $ is not expansive on $
X_{\lambda} $ i.e. there are $ n_{k} \rightarrow \infty, $ $ x_{k} \in X_{\lambda},
$ such that $ \vert (f^{n_{k}})^\prime (x_{k})\vert \leq 1. $ Now   using the
compactness of $ X_{\lambda} $ and the arguments above, we obtain a contradictions.
%%that any accumulation point of $\{f^{n_{k}}(x_{k})\}_{k \in {\N}}$ is
%%in the $w-limit$ $set$ of some recurrent critical point. Under our hypothesis that
%%proves that $f_{\lambda}$ is expansive.
Expansivity immediately implies summability of the point $ z = 1 $ and completes the
theorem.

\subsection{Proof of proposition 1}
%Does not exist a map $ f_{\lambda} $ such that
 %$ lim_{n \rightarrow \infty}|(f_{\lambda}^{n})^\prime(1)| = C > 0. $
%\end{prop}

\begin{proof}
We have $lim_{n \rightarrow
\infty}|\frac{1}{f_{\lambda}^{n+1'}(1)}/\frac{1}{f_{\lambda}^{n'}(1)}|=1$. Since
$lim_{n \rightarrow \infty}|\frac{f_{\lambda}^{n'}(1)}{f_{\lambda}^{n+1'}(1)}| = $\\
$ = lim_{n \rightarrow \infty}|\frac{1}{{\lambda}f_{\lambda}^{n+1}(1)}|, $ then $
|f_{\lambda}^{n}(1)|$ is near $\frac{1}{|{\lambda}|} $ for all large values of $ n.
$

This implies that $ X_{\lambda} $ is bounded, hence compact and $0$ is
non-recurrent, by Theorem 3, $ f_{\lambda} $ is summable. That is a contradiction
with the hypothesis.
\end{proof}

\section{Ruelle Operator: Definitions and Properties}

For any $ \lambda \in {\C^*} $ we define the following operators (compare with
\cite{Mak}, \cite{Mak1}, \cite{Lev}).

\begin{defi}\hfill\newline

\begin{itemize}

\item {\bf Ruelle operator (or push-forward operator)} $$ R_\lambda^* (\varphi)(z)
:= \sum_{\xi_i} \varphi(\xi_i) {\xi_i^\prime}^2 = \frac{1}{\lambda^2 z^2}
\sum_{\xi_i} \varphi(\xi_i), $$ where the summation is taken over all branches
$\xi_i$ of $ f_{\lambda}^{-1}. $ \item {\bf Modulus of the Ruelle operator}  $ \vert
R_\lambda^*\vert (\varphi)(z) = \frac{1}{\vert \lambda^2 z^2 \vert}\sum_{\xi_i}
\varphi(\xi_i). $ \item {\bf Beltrami operator} $ B_\lambda (\varphi) = \varphi
(f_\lambda) \frac{\overline{f_{\lambda}^\prime}}{f_{\lambda}^\prime}.$
\end{itemize}
\end{defi}

Then we have the following simple lemma.

\begin{lema}
For all $\lambda$ ;
\begin{enumerate}
\item $ R_\lambda^* : L_1 ({\C}) \rightarrow L_1({\C}) $ and $ \parallel R_\lambda^*
\parallel_{L_1} \leq 1, $
\item $ \vert R_\lambda^*\vert: L_1 ({\C}) \rightarrow L_1({\C}), $ $ \parallel
\vert R_\lambda^*\vert
\parallel_{L_1} \leq 1 $ and the fixed points of $ \vert R^*_{\lambda}\vert $ define
a finite, complex-valued, invariant, and absolutely continuous measures on $ {\C}.$

 \item $ B_\lambda : L_\infty ({\C}) \rightarrow L_\infty ({\C}), $ is the
dual operator to $ R_\lambda^*, $ and $ \parallel B_\lambda\parallel_{L_\infty} = 1.
$
\end{enumerate}
\end{lema}

\begin{proof} Immediately follows from the definitions.
\end{proof}

\subsection{Potential of Deformations}

The open unit ball $B$ of the space $Fix(B_{\lambda}) \subset L_{\infty}({\C})$ of
fixed points of $B_{\lambda}$ is called the space of invariant Beltrami
differentials for $f_{\lambda}$ and describes all quasiconformal deformations of
$f_{\lambda}$.

For $\mu \in B$ and for any $t$ with $ \vert t \vert < \frac{1}{\parallel \mu
\parallel}$, the element $\mu_{t}=t \mu \in B$. Let us denote by $h_{t}$ their
corresponding quasiconformal maps; then we have the following functional equation as
explained in \cite{Mak}, \cite{Mak1}:
$$
F_{\mu}(f_{\lambda}(z))-f_{\lambda}'(z)F_{\mu}(z)=G_{\mu}(z)
$$
where $h_{t}\circ f_{\lambda}\circ h_{t}^{-1}=f_{\lambda(t)} \in M_{f_{1}} $ and $
G_{\mu}(z) = \frac{\partial f_{\lambda(t)}}{\partial t}(z)\vert_{t=0} =
z\exp(\lambda z)\lambda^\prime(t)\vert_{t=0}.$ The function
$$
F_{\mu}(a)=\frac{\partial h_{t}}{\partial {t}}\vert_ {t=0}=-\frac{a(a-1)}{\pi}
\iint_{\C}\frac{\mu(z)}{z(z-1)(z-a)}
$$
is called the {\sl potential of the qc-deformations} generated by $ \mu $ and $
\overline{\partial}F_\mu = \mu $ in the sense of distributions, see \cite{G}.

\begin{lema}
If $ F(f_\lambda) = \emptyset, $ then $ G_\mu = 0 $ if and only if $\mu = 0. $
\end{lema}
\begin{proof}
If $ G_\mu = 0, $ then $ F_\mu(f_{\lambda}(z)) = f_{\lambda}'(z)F_{\mu}(z). $ Hence
$ F_\mu = 0 $ on the set of repelling periodic points and hence $ F_\mu = 0 $ on the
Julia set. Then  $ \mu = \overline{\partial}F_\mu = 0. $ The lemma is finished.
\end{proof}

Then by an inductive argument we have that

$$
 F_{\mu}(f^{n}_{\lambda}(a))=f_{\lambda}^{n'}(a) \left(
F_{\mu}(a)+\sum_{i=1}^{n}\frac{G_{\mu}(f_{\lambda}^{i-1}(a))} {f_{\lambda}^{i'}(a)}
\right ).
$$
from above $ G_{\mu}(a)=\frac{a f'_{\lambda}(a)c}{\lambda}, $ where the constant $ c
= \lambda^\prime(t)\vert_{t=0} $ and by the lemma 4 above $ c\neq 0. $

\begin{equation}
F_{\mu}(f^{n}_{\lambda}(a)) = f_{\lambda}^{n'}(a) \left (F_{\mu}(a) + \frac{a
c}{\lambda} + \frac{c}{\lambda^2} \sum_{i = 2}^{n}\frac{1}{(f_{\lambda}^{i -
2})^\prime(a)} \right) \label{ecu0}
\end{equation}

%If we redefine $F$ as $Q_{\mu}(a)=F_{\mu}(a)/F_{\mu}(f_{\lambda}(1))$, we %have a normalized potential

%with $Q_{\mu}(f_{\lambda}(1))=1$, with the above properties. So from now on %we can

%assume that the potential is normalized and will be denoted by $F_{\mu}$ %also.

%Let us note that in the limit as $n$ $\rightarrow$  $\infty$, %$f_{\lambda}^{n'}(a)$ $\rightarrow$

%$\infty$ if $a \in J(f)$, so %$F_{\mu}(a)+\frac{a}{f'_{\lambda}(1)}+\frac{1}{\lambda

%f'_{\lambda}(1)}(\sum_{i=2}^{n}\frac{1}{f_{\lambda}^{i-2'}(a)})$ must tend %to zero.

%{\bf Note:} A calculation of the last equation for $a=f_{\lambda}(1)$, %implies that $F_{\mu}(f_{\lambda}(1)+\frac{1}{\lambda}

%(1+\frac{1}{f'_{\lambda}(1)}+...+\frac{1}{f^{n'}_{\lambda}(1)})$ tends to %zero, since $F_{\mu}(f_{\lambda}(1))=1$ we have that

%$1+\frac{1}{f'_{\lambda}(1)}+...+\frac{1}{f^{n'}_{\lambda}(1)}$ tends to $-%{\lambda}$

%Next we have an application of the above concepts:

%%Now consider

%%$$
%%S_{n}=1 + \frac{1}{\lambda}\sum_{i = 2}^{n}\frac{1}{(f_{\lambda}^{i - 2})^\prime(1)}
%%$$

Now we are ready to prove the theorem 1. \subsection{Proof of Theorem 1} Firstly we
show (3). Such that $\frac{\vert (f^{n + 1}_{\lambda})^\prime(0)\vert}{\vert
(f^{n}_{\lambda})^\prime(0)\vert} = \vert\lambda f^{n + 1}_{\lambda}(0)\vert, $ then
assumption either
$$
\lim_{n\to \infty} sup \frac{\vert (f^{n + 1}_{\lambda})^\prime(1)\vert}{\vert
(f^{n}_{\lambda})^\prime(1)\vert} \leq C <\infty
$$
or
$$
\lim_{n\to \infty}inf\frac{\vert (f^{n +1}_{\lambda})^\prime(1)\vert}{\vert
(f^n_{\lambda})^\prime(1)\vert} > 0.
$$
implies either $ X_{\lambda} $ is a compact subset of the plane  or $ 0 \notin
X_{\lambda}, $ respectively. Assume $ F(f_{\lambda}) = \emptyset, $ then an
application of the theorem 3 implies a contradiction with $ \lim_{n\to \infty}\vert
(f^n_{\lambda})^\prime(0)\vert = 0. $ Hence we are done.

Now we show (1) and (2). Assume $f_{\lambda}$ is stable.

From the equation (1) above, we have that

$$
\frac{F_{\mu}(f^{n}_{\lambda}(a))}{(f_{\lambda}^{n})^\prime (a)}= F_{\mu}(a)+\frac{a
c}{\lambda}+\frac{c}{\lambda^2} \sum_{i=2}^{n}\frac{1}{(f_{\lambda}^{i - 2})^\prime
(a)}
$$

From \cite{G} we have the following inequality

$$
\vert F_{\mu}(a)\vert \leq M \vert a\vert \vert\log\vert a\vert\vert,
$$
where M is a constant depending only on $ \mu. $ Applying this estimate above we
obtain:

$$
\frac{\vert F_{\mu}(f^{n}_{\lambda}(a))\vert}{\vert (f_{\lambda}^{n})^\prime
(a)\vert}\leq \frac{M \vert f_{\lambda}^{n}(a)\vert\vert\log\vert
f^{n}_{\lambda}(a)\vert\vert} {\vert (f^{n}_{\lambda})^\prime (a)\vert}.
$$
Easy calculation shows $ log\vert f_{\lambda}^{n}(a)\vert = \vert
{\lambda}f_{\lambda}^{n-1}(a)\vert $ and \\
$ (f_{\lambda}^{n})^\prime (a) = {\lambda}^{2}f_{\lambda}^{n}(a)
f_{\lambda}^{n-1}(a)(f_{\lambda}^{n - 2})^\prime (a). $ Hence $ \frac{\vert
F_{\mu}(f^{n}_{\lambda}(a))\vert}{\vert (f_{\lambda}^{n})^\prime (a)\vert}\leq
\frac{M}{{\lambda}(f_{\lambda}^{n - 2})^\prime (a)}. $

Now let $ n_j $ be the sequence from the assumptions of theorem 1 items (1)-(2) and
the point $ a = 1. $ Since $ F_{\mu}(1) = 0, $ then from the equation 2 we obtain
the following equation:
$$
\frac{F_\mu(f_{\lambda}^{n_j + 2}(1))}{(f_{\lambda}^{n_j + 2})^\prime(1)} =
\frac{c}{\lambda}+ \frac{c}{\lambda^2} \sum_{i=2}^{n_j +
2}\frac{1}{(f_{\lambda}^{i-2})^\prime(1)} = \frac{c}{\lambda}\cdot S_{n_j}.
$$
Then this equation produces a contradiction in the both cases with the hypothesis
over $S_{n_j}$, so $f_{\lambda}$ is unstable.

\section{Calculation of the Ruelle Operator}

 In this section we calculate the action of the Ruelle operator on the family of rational
 functions  $\gamma_{a} (z) = \frac{a(a-1)}{z(z-1)(z-a)}$, such that $a \neq 0, 1$. Let us
 recall that any rational
integrable differential is a linear combination of such $\gamma_{a} (z)$.

Let $S = {\C} \backslash \{0,1\}$ be the trice punctured sphere.

\begin{prop}
\[ R_{\lambda}^{*} (\gamma_{a} (z)) = \frac{1}{(f_{\lambda})^\prime (a)}
\gamma_{f_{\lambda}(a)} (z) - \frac{a}{(f_{\lambda })^\prime (1)}
\gamma_{f_{\lambda}(1)} (z).
\]
\end{prop}

\begin{proof} Let $ h_a(z) = R_\lambda^* (\gamma_a)(z) - \frac{1}{(f_{\lambda})^\prime (a)}
\gamma_{f_{\lambda}(a)} (z) + \frac{a}{(f_{\lambda})^\prime (1)}
\gamma_{f_{\lambda}(1)} (z) $ be a function. Our aim is to show that $ h_a(z) $
defines a holomorphic integrable function  on the surface $ S, $ hence $ h_a(z) = 0
$ and we are done. By the lemma 3 the function $ h_a(z) $ is integrable over the
plane. Therefore it is enough to show that $ h_a(z) $ is holomorphic on  $ S. $

Let $ \varphi \in C^\infty(S) $ be any differentiable function with compact support
in $ S. $ Then
$$ \iint_{\C} \varphi_{\overline{z}} h_a(z) =
\iint_{\C} B_\lambda (\varphi_{\overline{z}}) \gamma_a(z) - \frac{1}{(f_{\lambda
})^\prime (a)}\iint_{\C}\varphi_{\overline{z}} \gamma_{f_{\lambda}(a)} (z) +
\frac{a}{(f_{\lambda})^\prime (1)}\iint_{\C} \varphi_{\overline{z}}
 \gamma_{f_{\lambda}(1)} =
$$

$$
= \iint_{\C} \varphi_{\overline{z}} (f_\lambda)
\frac{\overline{(f_{\lambda})^\prime}}{(f_{\lambda})^\prime} \gamma_a(z)-
\frac{1}{(f_{\lambda})^\prime (a)}\iint_{\C}\varphi_{\overline{z}}
\gamma_{f_{\lambda}(a)} (z) + \frac{a}{(f_{\lambda})^\prime (1)}\iint_{\C}
\varphi_{\overline{z}}
 \gamma_{f_{\lambda}(1)} = (\ast)
$$
On the other hand
$$
\iint_{\C} \varphi_{\overline{z}} (f_\lambda)
\frac{\overline{(f_\lambda)^\prime}}{(f_\lambda)^\prime} \gamma_a(z) =
 a(a -1) \iint_{\C}\frac{(\varphi \circ f_{\lambda})_{\overline{z}}} {z(z - 1)(z -
a)(f_{\lambda})^\prime} =(a - 1)\iint_{\C} \frac{(\varphi \circ
f_\lambda)_{\overline{z}}}{z (f_\lambda)^\prime}
$$
$$
- a \iint_{\C} \frac{(\varphi \circ f_\lambda)_{\overline{z}}} {(z - 1)
(f_\lambda)^\prime} + \iint_{\C} \frac{(\varphi \circ f_\lambda)_{\overline{z}}} {(z
- a) (f_\lambda)^\prime}.
$$

Such that always $\varphi (0) = 0, \; \varphi (1) = 0.$  Then

\[ (a - 1) \iint_{\C} \frac{(\varphi \circ f_\lambda)_{\overline{z}}}{z
(f_\lambda)^\prime} = \frac{a - 1}{f_\lambda^\prime(0)} \varphi(f_\lambda(0)) = 0
\]

\[ a\iint_{\C} \frac{(\varphi \circ f_\lambda)_{\overline{z}}}{(z-1)
(f_\lambda)^\prime} = \frac{a}{(f_\lambda)^\prime (1)} \varphi(f_\lambda(1))
\]

\[\iint_{\C} \frac{(\varphi \circ f_\lambda)_{\overline{z}}}{(z-a)
(f_\lambda)^\prime} = \frac{1}{(f_\lambda)^\prime (a)} \varphi(f_{\lambda} (a))
\]
the same decompositions show
\[\frac{1}{(f_{\lambda
})^\prime (a)}\iint_{\C}\varphi_{\overline{z}} \gamma_{f_{\lambda}(a)}(z) =
\frac{1}{(f_{\lambda})^\prime (a)} \varphi(f_\lambda(a))
\]

\[\frac{a}{(f_{\lambda})^\prime (1)}\iint_{\C} \varphi_{\overline{z}}
 \gamma_{f_{\lambda}(1)} = \frac{a}{(f_\lambda)^\prime(1)} \varphi(f_\lambda(1))
\]
and as a result we obtain
$$
(\ast) = 0.
$$
By the Weyl's Lemma $ h_a(z) $ is a holomorphic function on $ S. $ Hence we are
done.

\end{proof}

\begin{corol} If $ F(f_\lambda) =\emptyset $ and $ \mu\neq 0 \in B, $ then
$$
G_{\mu}(a) = \frac{af_{\lambda}^\prime(a)}{f_{\lambda}^\prime (1)} F_{\mu}
(f_{\lambda}(1)) $$
%%We can now calculate $)$ as defined in section 2.
\end{corol}
\begin{proof}
Let $ \mu\neq 0 \in B $ be invariant Beltrami differential for $ f_{\lambda}, $ then
by the proposition 3 we have

\[-\pi F_\mu(a) =  \iint \gamma_{a}(z) \mu = \iint R^*_{\lambda}(\gamma_{a}(z)) \mu  =
\frac{1}{f_{\lambda}^\prime (a)} (-\pi) F_{\mu} (f_{\lambda}(a)) -
\frac{a}{f_{\lambda}^\prime (1)} (-\pi) F_{\mu} (f_{\lambda}(1)). \]

Hence
$$ F_{\mu}(a) = \frac{1}{f_{\lambda}^\prime(a)} F_{\mu} (f_{\lambda}(a)) - \frac{a}
{f_{\lambda}^\prime (1)} F_{\mu} (f_{\lambda}(1)), $$ and

%\begin{equation}

$$ G_{\mu}(a)=F_{\mu}(f_{\lambda}(a))-f_{\lambda}^\prime(a) F_{\mu}(a) =
\frac{af_{\lambda}^\prime(a)}{f_{\lambda}^\prime (1)} F_{\mu} (f_{\lambda}(1)). $$

%\end{equation}
\end{proof}

%%We can choose $F_{\mu} (f_{\lambda}(1))=1$, hence the result is the equation

%%\begin{equation}
%%G_{\mu}(a)= \frac{af_{\lambda}^{'}(a)} {f_{\lambda}^{'} (1)}. \label{ecu3}
%%\end{equation}

Now, by the linearity of the Ruelle operator together with an easy induction argument, for any $n \geq 0$ we have

$$ (R_{\lambda}^{*})^{n} (\gamma_{a}(z)) = \frac{1}{(f^{n}_{\lambda })^\prime(a)}
\gamma_{f_{\lambda}^{n}(a)}(z) - \frac{f^{n - 1}_{\lambda }(a)}{(f^{n - 1}_{\lambda
})^\prime (a) f_{\lambda}^\prime(1)} \gamma_{f_{\lambda}(1)}(z)  - \leqno(*) $$

$$ -\frac{f^{n - 2}_{\lambda }(a)}{(f^{n - 2}_{\lambda })^\prime (a) f_{\lambda}^\prime(1)}
R_{\lambda}^{*}( \gamma_{f_{\lambda}(1)}(z) ) - \ldots -
\frac{a}{f_{\lambda}^\prime(1)} (R_{\lambda}^{*})^{n-1}(\gamma_{f_{\lambda}(1)}(z))
.
$$

Define the following series

$$
B(a) = \frac{1}{f_{\lambda}^\prime(1)} \sum_{j = 1}^{\infty} \frac{f^{j -
1}_{\lambda }(a)} { (f^{j - 1}_{\lambda })^\prime (a)}.
$$

\section{proof of the theorem 2 }

Assume $ f_\lambda $ is a stable map, then the summability of the singular value
implies $ F(f_{\lambda}) =\emptyset. $

Let  $ \mu \neq 0 $ be an invariant Beltrami differential. Then the formula $ (\ast)
$ above, the invariance of $ \mu, $ and
%%we have for any $a \in {\C}$
%%$$ \iint \gamma_{a}(z)\mu = \iint  R_{\lambda}^{*n}\gamma_{a}(z)\mu
%%$$
the definition of the potential $F_{\mu}$ give the following

$$F_{\mu}(a) =
\frac{1}{f_{\lambda}^{n'}(a)}F_{\mu}(f_{\lambda}^{n}(a))-B_n(a)F_{\mu}(f_{\lambda}(1)),\leqno(**)
$$
where $ B_n(a) $ is the $ n-th $ partial sum of the series $ B(a) $ above.

Let $ a $ be  a summable point, then the series $ B(a) $ is absolutely convergent
and by the arguments of the theorem 1, item (1), the expression
$\frac{1}{f_{\lambda}^{n'}(a)}F_{\mu}(f_{\lambda}^{n}(a)) \to 0$ as $n \to \infty. $

Then passing to the limit in the formula $ (\ast\ast) $ above  we have:

$$ F_{\mu}(a) = -B(a) F_{\mu} (f_{\lambda}(1)) $$

Now set $ a = f_{\lambda}(1), $ then:

\[ F_{\mu}(f_{\lambda}(1)) \left ( 1 + B ( f_{\lambda}(1)) \right ) = 0 . \]

and we have two possibilities:

1) $F_{\mu}(f_{\lambda}(1)) = 0$

Then by the Corollary 1, $ G_\mu = 0 $ and by the lemma 4, $ \mu = 0 $ which
contradicts the assumption above.

 2)$ B ( f_{\lambda}(1))= -1. $

Now we finish the theorem 2 in 3 steps. Let $ \varphi $ be the following series
$$
\varphi (z) := \sum_{n \geq 0} \frac{1}{(f_{\lambda}^{n})^\prime(f_{\lambda} (1))}
\gamma_{f_{\lambda}^{n}(f_{\lambda}(1))}(z), $$ then summability of the point $ z =
0 $ implies $ \varphi \in L_1({\C}). $

In the first step we show that under assumption 2) above, the function $
\vert\varphi\vert $ presents a density of a finite, invariant measure which is
absolutely continuous with respect to Lebesgue measure on the plane.
%%Then
%%we have the following result:

\begin{lema}Under assumption (2) above we have:
%%The function $ \varphi $ is a fixed point for Ruelle operator, that is
\[ R_{\lambda }^{*}(\varphi (z)) = \varphi (z) .\]
\end{lema}

\begin{proof}

%%\[ R_{\lambda }^{*} (\gamma_{f_{\lambda}(1)}(z)) = \frac{1}{f_{\lambda}^{'}(f_{\lambda}(1))}
%%\gamma_{f_{\lambda}(f_{\lambda}(1))}(z) - \frac{f_{\lambda}(1)}{f_{\lambda}^{'}(1)}
%%\gamma_{f_{\lambda}(1)}(z) , \]

For any $ n \geq 0, $ by the formula $ (\ast\ast) $  we have the following
expression

\[ R_{\lambda }^{*} \left (\frac{1}{(f_{\lambda}^{n})^\prime(f_{\lambda}
(1))} \gamma_{f_{\lambda}^{n}(f_{\lambda}(1))}(z) \right ) =
\]
\[\frac{1}{(f_{\lambda}^{n +1})^\prime(f_{\lambda}(1))} \gamma_{f_{\lambda}^{n
+1}(f_{\lambda}(1))}(z) - \frac{1}{f_{\lambda}^\prime(1)} \gamma_{f_{\lambda}(1)}(z)
\frac{f_{\lambda}^{n}(f_{\lambda }(1))}{(f_{\lambda}^{n})^\prime(f_{\lambda }(1))} ,
\]

Then summation over  all $n \geq 0$ gives

\[ R_{\lambda }^{*}(\varphi) = R_{\lambda }^{*} \left ( \sum_{n \geq 0}
\frac{\gamma_{f_{\lambda}^{n}(f_{\lambda}(1))}(z)}{(f_{\lambda}^{n})^\prime(f_{\lambda}
(1))} \right ) =  \]
\[ = \sum_{n \geq 0} \frac{1}{(f_{\lambda}^{n + 1})^\prime(f_{\lambda}
(1))} \gamma_{f_{\lambda}^{n + 1}(f_{\lambda}(1))}(z) -
\frac{1}{f_{\lambda}^\prime(1)} \gamma_{f_{\lambda}(1)}(z) \sum_{n \geq 0}
\frac{f_{\lambda}^{n}(f_{\lambda }(1))}{(f_{\lambda}^{n})^\prime(f_{\lambda }(1))} =
\]

\[ = \varphi (z) - \gamma_{f_{\lambda}(1)}(z) - \gamma_{f_{\lambda}(1)}(z) \left [
B(f_{\lambda}(1)) \right ]  = \varphi (z) \]

by hypothesis.

%%In conclusion $R_{\lambda }^{*}(\varphi (z)) = \varphi (z)$

\end{proof}

\begin{lema}In assumption of the lemma 5 above the function $|\varphi |$ is a fixed point
for the modulus of the Ruelle operator,
\[ | R_{\lambda }^{*}| (|\varphi |) = \vert \varphi \vert  . \]
\end{lema}

\begin{proof}

We recall that by definition, for every function $\varphi $

\[ | R_{\lambda }^{*}| (|\varphi |) = \sum_{\zeta_{i}} | \varphi (\zeta_{i})| \;
|\zeta_{i}^\prime |^{2} \]

\noindent where summation is over all branches $\zeta_{i}$ of inverses of
$f_{\lambda }(z) = e^{\lambda  z}$.

By assumption

$$
%%\begin{displaymath}
\parallel  \varphi  \parallel = \parallel R_{\lambda }^{*} (\varphi )
\parallel = \dobleint_{\C} \left | \sum_{\zeta_{i}} \varphi (\zeta_{i})
(\zeta_{i}^\prime)^{2} \right |.
%%\end{displaymath}
$$

Now define for each index $i$, $ \alpha_{i} = \varphi (\zeta_{i})
(\zeta_{i}^\prime)^{2}$, $\beta_{i} = \sum_{j \neq i} \varphi
(\zeta_{j})(\zeta_{j}^\prime)^{2}=\varphi-\alpha_{i}$.

With this notations we have

\[ \parallel \varphi \parallel = \parallel R_{\lambda }^{*} (\varphi )\parallel =
\dobleint_{\C} \left | \sum_{\zeta_{i}} \varphi (\zeta_{i}) (\zeta_{i}^\prime)^{2}
\right | dz \wedge d \bar z = \dobleint_{\C} |\alpha_{i} + \beta_{i}| \leq
\]

\[ \leq  \dobleint_{\C} | \alpha_{i} | +
\dobleint_{\C} | \beta_{i} |  \leq \dobleint_{\C} \sum_{\zeta_{i}}\vert \varphi
(\zeta_{i})(\zeta_{i}^\prime)^{2}\vert = \parallel \varphi
\parallel , \]

\noindent Hence all inequalities above are really equalities, then for each index
$i$ we have

\[ \dobleint_{\C} |\alpha_{i} + \beta_{i}|  =
\dobleint_{\C} | \alpha_{i} |  + \dobleint_{\C} | \beta_{i} |
\]

\noindent which implies that $|\alpha_{i} +  \beta_{i} | = |\alpha _{i}|+|\beta
_{i}|$ almost everywhere with respect to Lebesgue measure. Then for each index $i$

\[ |\alpha_{i} + \beta_{i}| = | \alpha_{i} + \sum_{j \neq i} \varphi (\zeta_{j})
(\zeta_{j}^\prime)^{2}| = |\alpha _{i}| + |\sum_{j \neq i} \alpha_{j} | , \] and by
the induction we obtain

$|\sum_{i} \alpha_{i}| = \sum_{i} |\alpha_{i} |. $

That implies that

\[ |\varphi | = |\sum_{i} \alpha_{i}| = \sum_{i} |\alpha_{i} | = \sum_{\zeta_{i}}
|\varphi (\zeta_{i})| |\zeta_{i}^\prime|^{2} = |R_{\lambda }^{*}|(|\varphi |). \]
\end{proof}

By the lemma 3 the measure  $\sigma(A) =  \iint_A \vert\varphi(z)\vert $ is a non -
negative invariant absolutely continuous probability measure, where $ A \subset
\chat $ is a measurable set. We have complete  the first step.

\medskip

Let $ Y={\C}-X_{\lambda} $ be the complement to the postsingular set $ X_{\lambda}.
$ In the second step we show that $ \varphi = 0 $ identically on $ Y. $

In the notation of the lemmas above we have:

\medskip

\begin{lema}
If $ \alpha_j \neq 0 $ identically on $ Y, $ then the function $ k_j =
\frac{\beta_j}{\alpha_j} $ is a non-negative constant on any component of  $  Y. $
\end{lema}

{\it Proof}. We have $ \vert 1 + \frac{\beta_j}{\alpha_j}\vert = 1 + {\biggl
|}\frac{\beta_j}{\alpha_j}{\biggr |}, $ then if $ \frac{\beta_j}{\alpha_j} =
\gamma^j_1 + i\gamma^j_2 $ we have
$$ {\biggl (}1 + (\gamma^j_1){\biggr )}^2 + (\gamma^j_2)^2 = {\biggl (}1 + \sqrt{(\gamma_1^j)^2 +(\gamma_2^j)^2}{\biggr )}^2 = 1 + (\gamma_1^j)^2 +(\gamma_2^j)^2 + 2\sqrt{(\gamma_1^j)^2 +(\gamma_2^j)^2}.
$$

Hence $ \gamma_2^j = 0  $ and $ \frac{\alpha_j}{\beta_j} = \gamma_1^j $ is a
real-valued function but $\frac{\alpha_j}{\beta_j}$ is meromorphic function. So $
\gamma_1^j = k_j $ is constant on every connected component of $ Y $ and the
condition $ \vert 1 + k_j\vert = 1 +\vert k_j\vert $ shows $ k_j \geq 0.$

\bigskip

\bigskip

\begin{defi}

A measurable set $A \in \chat$ is called back wandering if and only if $m(f^{-n}(A)
\cap f^{-k}(A)) = 0$, for $ k \neq n$.

\end{defi}

\begin{corol}

If $\varphi \neq 0$ on $Y$, then (i) $ m(X_{\lambda}) = 0, $ where $ m $ is the
Lebesgue measure and (ii) $ \frac{\bar{\varphi}}{\vert \varphi \vert} $  defines an
invariant Beltrami differential.

\end{corol}

\begin{proof}

(i) If $ m(X_{\lambda}) > 0, $  then $ m(f_{\lambda}^{-1}(X_{0})) > 0 $ so $
m(f_{\lambda}^{-1}(X_{0}) - X_{0}) > 0 $ since $ f_{\lambda}^{-1}(X_{0}) \neq X_{0},
$ $ X_{\lambda} \neq \C, $ denote by $ Z_1 =  f_{\lambda}^{-1}(X_{\lambda}) -
X_{\lambda}. $ Then $ Z_1 $ is back wandering thus $ \varphi = 0, $ on the orbit of
$ Z_1, $ which is dense in $ J(f_{\lambda}), $ hence $ \varphi = 0 $ in $ Y. $
Therefore, $ m(X_{\lambda}) = 0. $

\bigskip

(ii) By notations and the proofs of Lemmas 5 and 6 we have
  $k_i(x)= \frac{{\beta}_{i}}{{\alpha}_{i}}= \frac{{\varphi}}{{\alpha}_{i}}-1$
so ${\varphi(x)}=(1+k_{i}(x)){\alpha}_{i} =
(1+k_{i}(x))(\varphi(\zeta_i(x))(\zeta_i')^{2}(x) $. Hence,

$$
%%\begin{displaymath}
\frac{{\bar
{\varphi}(x)}}{|{\varphi}(x)|}=\frac{(1+k_{i}(x)){\bar{\varphi}({\zeta}_{i}(x))}{\bar{(\zeta_i')}^{2}(x)}}{(1+k_{i}(x))|{\varphi}({\zeta}_{i}(x)||(\zeta_i')^{2}(x)|},
%%\end{displaymath}
$$

\noindent and so for any branch $ \zeta_i $ we have

$$
%%\begin{displaymath}
{\mu}=\frac{{\bar {\varphi}}}{|{\varphi}|}=\frac{{\bar {\varphi}({\zeta}_{i})}{\bar
{\zeta}_{i}'}}{|{\varphi}({\zeta}_{i})|{\zeta}_{i}'} =\mu(\zeta_{i})\frac{\bar
{\zeta}_{i}'}{{\zeta}_{i}}
%%\end{displaymath}
$$

as result $\mu = \frac{\bar{\varphi}}{\vert\varphi\vert}$ is an invariant line
field. Thus the corollary is proved.

\end{proof}

\bigskip
Now we prove the main result of the second step.

\begin{prop}

If $\varphi \neq 0$ on $Y$, then $f_{\lambda}$ is unstable.

\end{prop}

\begin{proof}

Let us show first that $ X_{\lambda} = {\bigcup}f_{\lambda}^{i}(1). $ We will use a
McMullen argument as in \cite{MM}. By  Corollary 2, $
\mu=\frac{\bar{\varphi}}{|{\varphi}|} $ is an invariant Beltrami differential. That
implies that $ {\varphi} $ is dual to $ {\mu} $ and $ \varphi $ is defined by $ \mu
$ up to a constant. We will construct a meromorphic function $ \psi, $ dual to $ \mu
$ and such that $ \psi $ has finite number of poles on each disc $ D_{R} $ of radius
$ R $ centered at $ 0. $

For that suppose that for $ z \in \C $ there exists a branch $ g $ of a suitable $
f_{\lambda}^{n}, $ such that $ g(U_{z}) \in Y, $ where $ U_{z} $ is a neighborhood
of $ z. $ Then define $ \psi(\zeta)=\varphi(g(\zeta))(g')^{2}(\zeta), $ for all $
\zeta \in U_{z}. $ Note that $ \psi(\zeta) $ is dual to $ \mu $ and has no poles in
$ U_{z}. $ If there is no such branch $ g, $ then $ \zeta $ is in the postsingular
set, and there is a branched covering $ F $ from a neighborhood of $ \zeta $ to $
U_{z}, $ then define $ \psi(\zeta) = F^{*}(\varphi), $ with $ F^{*} $ the Ruelle
operator of $ F. $ The map $ \psi $ is a meromorphic function dual to $ \mu $ in $
U_{z} $ and has finite number of poles.

By considering $ R\to \infty $ we construct a  meromorphic function $ \psi $ which
is dual to $ \mu. $ The poles of $ \psi $ forms a discrete set accumulating to $ z =
\infty. $ Since  $ \varphi $ is a dual to $ \mu, $ then $ \varphi = C\cdot\psi, $
where $ C $ is a constant. Hence $ X_{\lambda} = {\bigcup}f^{i}(1) $ is a discrete
closed set accumulating to $ z = \infty $ and $ Y $ is connected.

By the Corollary 2 the functions $ k_{i} $ are globally defined constants on $ Y. $
Moreover by the argument of the lemma 7
$\varphi(x)=(1+k_{i})(\varphi(\zeta_i(x))(\zeta_i')^{2}(x) $ for any $ x \in {\C}, $
thus  $ k_i = k_j $ for any $ i, j. $

 So we have
$\sum_{i}\frac{\varphi(x)}{1+k_{i}}=\sum_{i}\varphi(\zeta_i(x))(\zeta_i')^{2}(x)=\varphi(x)$,
since the first term of the equation is infinite, this can be only iff $\varphi=0$.

\end{proof}

Now to obtain a contradiction, in the step 3 we show that if $ f_{\lambda} \in W $
is a structurally stable, then $ \varphi \neq 0 $ identically on $ Y. $

The following proposition is proved in \cite{Mak1}.

\begin{prop} Let $ a_i \in \C, a_i\neq a_j, $ for $ i\neq j $ be points such that
$ Z = \overline{\cup_i a_i} \subset {\C} $ is a compact set. Let  $ b_i \neq 0 $ be
complex numbers such that the series $ \sum b_i $ is absolutely convergent. Then the
function $ l(z) = \sum_i\frac{b_i}{z - a_i} \neq 0 $ identically on $ Y = {\C}
\backslash Z $ in any of the following cases

\begin{enumerate}
\item the set $ Z $ has zero Lebesgue measure \item if diameters of components of
$\C \backslash Z $ uniformly bounded below from zero and \item If $ O_j $ denote the
components of $ Y, $ then $ Z \in \cup_j{\partial O_j}.$
\end{enumerate}

\end{prop}

\begin{prop} Let  $ f_{\lambda}  $ be the exponential map and $ 0 $  a summable point.
 Then $\varphi(z) \neq 0 $ identically on $ Y$ in any of the following cases
\begin{enumerate}
\item if $ 0 \notin X_{\lambda}, $ \item if diameters of components of $ Y $ are
uniformly bounded below from 0, \item  If $ m\left (X_{\lambda}\right ) = 0, $ where
$ m $ is the Lebesgue measure on $ {\C}, $
%\item if $ X_{0} \subset \cup_i \partial D_i, $ where $ D_i $ are components of Fatou set.

\end{enumerate}

\end{prop}

{\it Proof.} Let us prove (1). Denote $d_{\lambda}=f_{\lambda}(1)$

Assume now that the set $ X_{\lambda} $ is bounded. Then by Proposition 4 we have
that $ \varphi(z) = \frac{C_1}{z}+\frac{C_2}{z-1}+ \sum\frac{1}{(f_{\lambda}
^i)'(d_{\lambda})(z - f_{\lambda} ^i(d_{\lambda}))} = l(z) \neq 0$. The other cases
follows directly from Proposition 4 also.

\par Now let $ X_{\lambda} $ be unbounded. Let $ y \in {\C} $ be a point such that the point
$ 1 - y \in Y, $ then the map $ g(z) = \frac{y z}{z + y - 1} $ maps $ X_{\lambda} $
into $ \C. $ Let us consider the function $ G(z) =
\frac{1}{z}\sum_i\frac{(f_{\lambda} ^i(d_{\lambda}) - 1)}{(f_{\lambda}
^i)'(d_{\lambda})} - \frac{1}{z - 1}\sum_i\frac{f_{\lambda}
^i(d_{\lambda})}{(f_{\lambda} ^i)'(d_{\lambda})} + \sum\frac{1}{(f_{\lambda}
^i)'(d_{\lambda})(z - g(f_{\lambda} ^i(d_{\lambda}))}, $ then by proposition 4 $
G(z) \neq 0 $ identically on $ g(Y). $
\bigskip

 Now we {\bf Claim} that  $G(g(z))g'(z) = \phi(z).$

\par {\it Proof of the claim.} Let us define $ C_1 = \sum_i\frac{(f_{\lambda} ^i(d_{\lambda}) - 1)}{(f_{\lambda} ^i)'(d_{\lambda})} $ and $ C_2 = \sum_i\frac{f_{\lambda} ^i(d_{\lambda})}{(f_{\lambda} ^i)'(d_{\lambda})} $ then we have

$$
\frac{C_1}{g(z)} = \frac{C_1(z + y - 1)}{yz} \text{ and } \frac{C_2}{g(z) - 1} =
\frac{C_2(z + y - 1)}{(y - 1)(z - 1)}
$$

and for any $ n $

%%\begin{displaymath}
\begin{multline}
\frac{1}{g(z) - g(f_{\lambda} ^n(d_{\lambda}))} = \frac{(z + y - 1)(f_{\lambda}
^n(d_{\lambda}) + y - 1)}{y(y - 1)(z - f_{\lambda} ^n(d_{\lambda}))} =\nonumber\\
= \frac{1}{y(y - 1)}\left(\frac{(z + y - 1)^2}{z - f_{\lambda} ^n(d_{\lambda})} + 1
- y - z\right),
%%\end{displaymath}
\end{multline}

then

\begin{multline}
G(g(z)) = \frac{C_1(z + y - 1)}{yz} -  \frac{C_2(z + y - 1)}{(y - 1)(z - 1)} + \sum \frac{1}{(f_{\lambda} ^i)'(d_{\lambda})(g(z) - g(f_{\lambda} ^i(d_{\lambda}))} =\nonumber\\
= \frac{1}{y(y - 1)}\biggl((1 - y - z)\sum \frac{1}{(f_{\lambda} ^i)'(d_{\lambda})} + (z + y - 1)^2\sum\frac{1}{(f_{\lambda} ^i)'(d_{\lambda})(z - f_{\lambda} ^i(d_{\lambda}))} + \nonumber\\
+ \frac{C_1(z + y - 1)}{yz} -  \frac{C_2(z + y - 1)}{(y - 1)(z - 1)}\biggr) =
\nonumber\ast
\end{multline}

and

\begin{gather}
\ast = \frac{1}{g'(z)}\left( \phi(z) + \frac{\sum_i\frac{f_{\lambda} ^i(d_{\lambda})
- 1}{(f_{\lambda} ^i)'(d_{\lambda})}}{z} - \frac{\sum_i\frac{f_{\lambda}
^i(d_{\lambda})}{(f_{\lambda} ^i)'(d_{\lambda})}}{z - 1} +
\frac{\sum\frac{1}{(f_{\lambda} ^i)'(d_{\lambda})}}{1 - y - z}\right) + \nonumber\\
+  \frac{1}{g'(z)}\left( \frac{C_1(y - 1)}{z(z + y - 1)}
 - \frac{C_2y}{(z - 1)(z + y - 1)}\right )= \nonumber\\
= \frac{\phi(z)}{g'(z)}.\nonumber
\end{gather}

Hence $ \phi(z) = 0 $ identically on $ Y $ if and only if  $ G(z) = 0 $ identically
on $ g(Y).$ So by proposition 4 we complete the proof of this proposition.

Step 3 and the theorem 2 are finished.

%%\begin{theo} Let $f_{\lambda} \in X$. Then $f_{\lambda} $ is an unstable map hence does not have invariant Beltrami differentials in its Julia set.
%%\end{theo}

%%{\it Proof.} Assume $ f_{\lambda}  $ is an stable map.

%%By the lemma 4, there exist $\varphi $, such that $ R_{\lambda} ^*(\varphi) =
%%\varphi $ and $ \varphi \neq 0 $ identically on $ Y $ by the proposition 5. Then by
%%an application of the proposition 3 we obtain a contradiction with the assumption of
%%stability.  Now, an invariant Beltrami differential implies a one dimensional set of
%%deformations but that is impossible because the full family is one dimensional and
%%unstable.

\bigskip

%%[1] Devaney R., Fagella N., Jarque X.; Hyperbolic components of the Complex
%%Exponential Family. Fundamenta Mathematicae 174 (2002),193-215.

\noindent{\tt Petr M. Makienko}\\
{Permanent addresses:}\\
{ \small Instituto de Matematicas, UNAM}\\
{\small Av. de Universidad s/N., Col. Lomas de Chamilpa}\\
{\small Cuernavaca, Morelos, C.P. 62210, M\'exico.}\\
{\small \tt E-mail: makienko@aluxe.matcuer.unam.mx}\\
 {\small and}\\
{\small Institute for Applied Mathematics,}\\
{\small   9 Shevchenko str.,}\\
{\small Khabarovsk, Russia}\\
{\small \tt E-mail makienko@iam.khv.ru}

\bigskip

\noindent{\tt Guillermo Sienra}\\
{ \small Facultad de Ciencias, UNAM}\\
{ \small Av. Universidad 30, C.U.}\\
{\small M\'exico D.F., C.P. 04510,   M\'exico.}\\
{\small \tt E-mail gsl@hp.fciencias.unam.mx}


\bigskip



\begin{thebibliography}{99}
\bibitem[1]{Av} A. Avila,   {\textit Infinitesimal perturbations of rational maps.}
Nonlinearity, N°15, p. 695-704, (2002).
\bibitem[2]{D} R. Devaney,  {\textit Structural Instability of $Exp(z).$}
Proceedings of the American Mathematical Society,  94 (1985), 545-548.
\bibitem[3]{DMS} P. Dominguez, P. Makienko, G. Sienra,  {\textit Ruelle operator and transcendental
entire maps.}  Discrete and Continuous Dynamics. Vol 12, N°4, p.773-789, (2005).
\bibitem[4]{DG}   A. Douady, L. R. Goldberg,  {\textit The nonconjugacy of certain exponential
functions.} Holomorphic functions and moduli, Vol. I (Berkeley, CA, 1986), Math.
Sci. Res. Inst. Publ., vol. 10, Springer, New York,  pp. 1–7, (1988).
\bibitem[5]{EL} A. Eremenko, M. Lyubich,  {\textit Dynamical properties of some classes of entire
functions.} Ann. Inst. Fourier, Grenoble. 42, p. 989-1020, (1992).
\bibitem[6]{Lev}    G. Levin,  {\textit On Analytic Approach to The Fatou
Conjecture} Fundamenta Mathematicae, 171, p. 177-196, (2002).
\bibitem[7]{Mak} P. Makienko, {\textit Remarks on Ruelle operator and line field
problem. I} Preprint. FIM, Zurich (2000).
\bibitem[8]{Mak1} P. Makienko, {\textit Remarks on Ruelle operator and line field problem II.} To be
published in Ergodic Theory.
\bibitem[9]{MM} C. McMullen,  {\textit Complex Dynamics and Renormalization.} Annals of Mathematical
Studies, Princeton University Press, (1994).
\bibitem[10]{MSS}  R. Man{\'e}, P. Sad and D.
Sullivan,  {\textit On the dynamic of rational maps.}  Ann. Sci. Ec. Norm. Sup.  16,
p. 193 - 217, (1983).
\bibitem[11]{G} F. P. Gardiner,  {\textit Teichmuller theory and quadratic differentials.} John Wiley and
Sons. NY. (1987).
\bibitem[12]{ST} M. Shishikura, T. Lei, {\textit An alternative proof of Ma\~{n}e's theorem on
non-expanding Julia sets.}  The Mandelbrot set, theme and variations, p. 265 - 279,
London Math. Soc. Lecture Note Ser., 274, Cambridge Univ. Press, Cambridge, (2000).
\bibitem[13]{Ye}    Zhuan Ye, {\textit Structural instability of exponential functions.} Trans. Amer. Math.
Soc. 344, no. 1, p. 379–389, (1994).
\end{thebibliography}
\end{document}